\documentclass[twoside,11pt,reqno]{amsart}

\usepackage{amsmath,amsthm,amssymb,amstext,amsfonts,amscd}
\usepackage{graphicx}
\usepackage{xcolor}
\usepackage{multirow}

\newtheorem{theorem}{Theorem}
\theoremstyle{plain}

\newtheorem{corollary}{Corollary}

\newtheorem{proposition}{Proposition}

\numberwithin{equation}{section}

\begin{document}
\title[\hfil Fractional differential and $k$-Mittag-Leffler function]
{On the solutions of certain fractional kinetic equations involving $E^{\gamma,q}_{k,\alpha,\beta}(.)$ }

\author[ P. Agarwal$^{1, 2}$, D. O'Regan$^{3}$ and M. Chand$^{4}$]{ Praveen
Agarwal, Donal  O'Regan and Mehar Chand}

\bigskip
\address{P. Agarwal:
$^{1}$Department of Mathematics,
Anand International College of Engineering,
Jaipur303012, India\\$^{2}$Department of Mathematics,
University Putra Malaysia,
43400 UPM, Serdang, Selangor, Malaysia}
 \email{goyal.praveen2011@gmail.com}

\bigskip
\address{D.O'Regan:$^{3}$School of Mathematics, Statistics and Applied Mathematics, National University of Ireland, Galway, IRELAND}
\email{donal.oregan@nuigalway.ie}

\bigskip
\address{M. Chand: $^{4}$Department of Mathematics, Fateh College for Women, Bathinda 151103, India}
\email{mehar.jallandhra@gmail.com}




\begin{abstract}
We develop a new generalized form of the fractional kinetic equation
involving a generalized k-Bessel function. The generalized
$k$-Mittag-leffler function $E^{\gamma,q}_{k,\alpha,\beta}(.)$ is
discussed in terms of the solution of the fractional kinetic
equation in the present paper. The results obtained here are quite
general in nature and capable of yielding known and as well new
results.
 \end{abstract}

\maketitle

\section{Introduction and Preliminaries}

The $k$-Pochhammer symbol  was introduced in \cite{1} as follows:
\begin{equation}\label{Poch}
\aligned & (x)_{n,k}
  =x(x+k)(x+2k)...(x+(n-1)k),\endaligned
\end{equation}
\begin{equation}\label{Poch-1}
\aligned & (x)_{(n+r)q,k}
  =(x)_{rq,k}(x+qrk)_{nq,k},\endaligned
\end{equation}
where $x\in\mathbb{C}, k\in\mathbb{R}$ and $n\in\mathbb{N}$.

\begin{proposition}\label{P-1.1} Let $\gamma\in\mathbb{C}, k,s\in\mathbb{R}$. Then the following identity holds
\begin{equation}\label{P-1}\aligned & \Gamma_s(\gamma)
  =\left(\frac{s}{k}\right)^{\frac{\gamma}{s}-1}\Gamma_k\left(\frac{k\gamma}{s}\right),\endaligned\end{equation}
and in particular
\begin{equation}\label{Poch-2}
\aligned
\Gamma_k(\gamma)=k^{\frac{\gamma}{k}-1}\Gamma\left(\frac{\gamma}{k}\right).
\endaligned
\end{equation}
\end{proposition}

\begin{proposition}\label{P-1.2} Let $\gamma\in\mathbb{C}, k,s\in\mathbb{R}$ and $\gamma\in\mathbb{C}$. Then the following identity holds
\begin{equation}\label{P-2}\aligned & (\gamma)_{nq,s}=\left(\frac{s}{k}\right)^{nq}\left(\frac{k\gamma}{s}\right)_{nq},\endaligned
\end{equation}
and in particular
\begin{equation}\label{Poch-3}
\aligned&
(\gamma)_{nq,k}=(k)^{nq}\left(\frac{\gamma}{k}\right)_{nq}.
\endaligned
\end{equation}
\end{proposition}
For more details on the $k$-Pochhammer symbol, the $k$-special function
and the fractional Fourier transform we refer the reader to the
papers by Romero et. al. \cite{110,116}.

Let $k\in \mathbb{R}, \alpha, \beta, \gamma\in \mathbb{C};
\Re(\alpha)>0, \Re(\beta)>0, \Re{(\gamma)}>0$ and $ q\in
(0,1)\cup\mathbb{N}$. Then the generalized $k$-Mittag-Leffler
function, denoted by $ E^{\gamma,q}_{k,\alpha,\beta}(z)$, is defined
as:
\begin{equation}\label{k-ML}
\aligned&
E^{\gamma,q}_{k,\alpha,\beta}(z)=\sum^\infty_{n=0}\frac{(\gamma)_{nq,k}
z^n}{\Gamma_k(n \alpha+\beta)n!}\,\,,
\endaligned
\end{equation}
where $(\gamma)_{nq,k}$ denotes the $k$-Pochhammer symbol given by
equation \eqref{Poch-3} and $\Gamma_k(x)$ is the $k$-gamma function
given by the equation \eqref{Poch-2} (also see \cite{118}).

The generalized Pochhammer symbol is defined as (see \cite[page 22]{Rainville-ED-1963})
\begin{equation}\label{G-poch}
\aligned& (\gamma)_{nq}=\frac{\Gamma(\gamma+nq)}{\Gamma(\gamma)}=q^{qn}\prod^q_{r=1}\left(\frac{\gamma+r+1}{q}\right)_n,~~~ {\rm if}~ q\in \mathbb{N}.
\endaligned
\end{equation}
We consider particular cases of $E^{\gamma,q}_{k,\alpha,\beta}(z)$:
\newline

\noindent (i) For $q=1$, equation\eqref{k-ML} yields the
$k$-Mittag-Leffler function (see \cite{2}), defined as:
\begin{equation}\label{k-ML-1}
\aligned&
E^{\gamma,1}_{k,\alpha,\beta}(z)=\sum^\infty_{n=0}\frac{(\gamma)_{n,k}
z^n}{\Gamma_k(n \alpha+\beta)n!}=E^{\gamma}_{k,\alpha,\beta}(z),
\endaligned
\end{equation}

\noindent (ii) For $k=1$, n\eqref{k-ML} yields the
Mittag-Leffler function, defined in \cite{24}:
\begin{equation}\label{k-ML-2}
\aligned & E^{\gamma,q}_{1,\alpha,\beta}(z)=\sum^\infty_{n=0}\frac{(\gamma)_{nq} z^n}{\Gamma(n \alpha+\beta)n!}=E^{\gamma,q}_{\alpha,\beta}(z),
\endaligned
\end{equation}

\noindent (iii) For $q=1$ and $k=1$, \eqref{k-ML} gives the
Mittag-Leffler function, defined in Dorrego and Cerutti \cite{2}:
\begin{equation}\label{k-ML-3}
\aligned&
E^{\gamma,1}_{1,\alpha,\beta}(z)=\sum^\infty_{n=0}\frac{(\gamma)_{n}
z^n}{\Gamma(n \alpha+\beta)n!}=E^{\gamma}_{\alpha,\beta}(z),
\endaligned
\end{equation}

\noindent (iv) For $q=1, k=1$ and $\gamma=1$, \eqref{k-ML}
gives the Mittag-Leffler function (see \cite{Wiman -1905}), defined
as:
\begin{equation}\label{k-ML-4}
\aligned& E^{1,1}_{1,\alpha,\beta}(z)=\sum^\infty_{n=0}\frac{
z^n}{\Gamma(n \alpha+\beta)n!}=E_{\alpha,\beta}(z),
\endaligned
\end{equation}

\noindent (v) For $q=1, k=1, \gamma=1$ and $\beta=1$,
\eqref{k-ML} gives the Mittag-Leffler function (see
\cite{20}), defined  as:
\begin{equation}\label{k-ML-5}
\aligned & E^{1,1}_{1,\alpha,1}(z)=\sum^\infty_{n=0}\frac{ z^n}{\Gamma(n \alpha+1)n!}=E_{\alpha}(z).
\endaligned
\end{equation}

A detailed account of the Mittag-Leffler function and their
applications can be found in the survey paper by Saxena et. al. \cite{25}.

\section{Fractional differential equations}
\noindent In physics, dynamical systems, control systems and
engineering, fractional differential equations model  many physical
phenomena and kinetic equations describe the continuity of motion of
substances. The extension and generalization of fractional kinetic
equations involving fractional operators can be found in
\cite{7, 17,  11}.

The fractional differential equation between the rate of change of
the reaction, the destruction rate and the production rate was
established by Haubold and Mathai\cite{7} and is given as follows:
\begin{eqnarray}\label{D-1} \frac{dN}{dt}=-d(N_{t})+p(N_{t}), \end{eqnarray}
where $ N=N(t) $ denotes the rate of reaction,  $ d=d(N) $ the rate
of destruction,   $ p=p(N) $ the rate of production and $N_{t} $
denotes the function defined by   $N_{t}(t^{*})= N(t-t^{*}),
t^{*}>0. $

The special case of \eqref{D-1} for spatial fluctuations and
inhomogeneities in $ N(t) $  the quantities are neglected , that is
the equation
\begin{eqnarray}\aligned &\label{D-2} \frac{dN}{dt}=-c_{i}N_{i}(t),\endaligned \end{eqnarray}
with the initial condition that $ N_{i}(t=0)=N_{0} $ is the number
density of the species $ i $ at time $ t=0 $ and $ c_{i}>0 $. If we
remove the index $ i $ and integrate the standard kinetic equation
\eqref{D-2}, we have
\begin{eqnarray}\aligned &\label{D-3} N(t)-N_{0}=-c{}_{0}D^{-1}_{t}N(t),\endaligned \end{eqnarray}
where $ {}_{0}D^{-1}_{t} $ is the special case of the
Riemann-Liouville integral operator $ {}_{0}D^{-\nu}_{t} $ defined
as
\begin{eqnarray}\aligned &\label{D-4} {}_{0}D^{-\nu}_{t}f(t)=\frac{1}{\Gamma(\nu)}\int_{0}^{t}\left(t-s\right)^{\nu-1}f(s)ds,\hskip 6mm  (t>0,R(\nu)>0).
\endaligned \end{eqnarray}
The fractional generalization of the standard kinetic
equation\eqref{D-3} is given by Haubold and Mathai \cite{7} as
follows:
\begin{eqnarray}\label{D-5} N(t)-N_{0}=-c^{\nu}{}_{0}D^{-1}_{t}N(t), \end{eqnarray}
and they obtained the solution of \eqref{D-5} as
\begin{eqnarray}\aligned &\label{D-6} N(t)=N_{0}\sum _{k=0}^{\infty }\frac{(-1)^{k}}{\Gamma\left(\nu k+1\right)}\left(ct\right)^{\nu k}.\endaligned \end{eqnarray}
In \cite{11} the authors considered the following fractional kinetic
equation:
\begin{eqnarray}\aligned &\label{D-7} N(t)-N_{0}f(t)=-c^{\nu}{}_{0}D^{-\nu}_{t}N(t),\hskip 6mm (\Re(v)>0),\endaligned\end{eqnarray}
where $ N(t) $ denotes the number density of a given species at time
$ t $, $ N_{0}=N(0) $ is the number density of that species at time
$t=0 $, $ c $ is a constant and $ f \in \mathcal{L}(0,\infty) $.
Apply the Laplace transform to \eqref{D-7} (see \cite{17}), so
\begin{eqnarray}\aligned &\label{D-8} L\left\{ N(t);p\right\}=N_{0}\frac{F(p)}{1+c^{\nu}p^{-\nu}}=N_{0}\left(\sum _{n=0}^{\infty }(-c^{\nu})^{n}p^{-\nu n}\right)F(p),\\&\hskip 6mm \left(n\in N_{0},\left|\frac{c}{p}\right|<1\right)\endaligned \end{eqnarray}
where the Laplace transform \cite{19} is given by
\begin{eqnarray}\aligned &\label{D-9} F(p)=L\left\{ N(t);p\right\} =\int_{0}^{\infty} e^{-pt}f(t)dt,\hskip 6mm (\mathcal{R}(p)>0).\endaligned  \end{eqnarray}

The objective of this paper is to derive the solution of the
fractional kinetic equation involving the generalized
$k$-Mittag-Leffler function. The results obtained in terms of the
Mittag-Leffler function are rather  general in nature and we can
easily construct various known and new fractional kinetic equations.

\section{Solution of generalized fractional kinetic equations}
In this section, we investigate the solution of the generalized
fractional kinetic equations by considering the generalized $k$-Bessel
function.

\begin{theorem}\label{Th-1} If $ d>0, \nu>0,$ $k\in \mathbb{R}, \alpha, \beta, \gamma\in \mathbb{C}; \Re(\alpha)>0, \Re(\beta)>0, \Re{(\gamma)}>0$ and $ q\in (0,1)\cup\mathbb{N}$, then the solution of the equation
\begin{eqnarray}\aligned &\label{Th-1-1} N(t)-N_{0}E^{\gamma,\tau}_{k,\alpha,\beta}(t)=-d^{\nu}{}_{0}D^{-\nu}_{t}N(t)\endaligned\end{eqnarray}
is given by
\begin{eqnarray}\aligned &\label{Th-1-2} N(t)=N_{0}\sum^\infty_{n=0}\frac{(\gamma)_{n\tau,k}}{\Gamma_k(n \alpha+\beta)}t^n E_{\nu,n+1}(-d^{\nu}t^{\nu}),
\endaligned \end{eqnarray}
where the generalized Mittag-Leffler function $ E_{\alpha,\beta}(x) $ is given by \cite{20}
\begin{eqnarray}\label{2.3}\aligned & E_{\alpha,\beta}(x)=\sum _{n=0}^{\infty }\frac{(x)^{n}}{\Gamma\left(\alpha n+\beta\right)}.\endaligned \end{eqnarray}

\end{theorem}

Proof: The Laplace transform of the Riemann-Liouville fractional
integral operator is given by \cite{21}
\begin{eqnarray}\aligned &\label{Th-1-3 } L\left\{{}_{0}D^{-\nu}_{t}f(t);p\right\}=p^{-\nu}F(p),\endaligned \end{eqnarray}
where $ F(p) $ is defined in \eqref{D-9}. Now, applying the Laplace
transform to both sides of \eqref{Th-1-1} gives
\begin{eqnarray}\aligned &\label{Th-1-4}L\left\{ N(t);p\right\}=N_{0}L\left\{ E^{\gamma,\tau}_{k,\alpha,\beta}(t);p\right\}-d^{\nu}L\left\{{}_{0}D^{-\nu}_{t}N(t);p\right\},\endaligned \end{eqnarray}
\begin{eqnarray}\aligned &\label{Th-1-5} N(p)=N_{0}\left(\int^{\infty}_{0}e^{-pt}\sum^\infty_{n=0}\frac{(\gamma)_{n\tau,k}}{\Gamma_k(n \alpha+\beta)}\frac{t^n}{n!}dt\right)\\&\hskip 6mm-d^{\nu}p^{-\nu}N(p),\endaligned\end{eqnarray}
\begin{eqnarray}\aligned &\label{Th-1-6} N(p)+d^{\nu}p^{-\nu}N(p)=N_{0}\sum^\infty_{n=0}\frac{(\gamma)_{n\tau,k}}{\Gamma_k(n \alpha+\beta)}\frac{1}{n!}\\&\hskip 6mm\times\int^{\infty}_{0}e^{-pt}t^{n}dt\endaligned \end{eqnarray}
\begin{eqnarray}\aligned &\label{Th-1-7} =N_{0}\sum^\infty_{n=0}\frac{(\gamma)_{n\tau,k}}{\Gamma_k(n \alpha+\beta)}\frac{1}{n!}\frac{\Gamma(n+1)}{p^{n+1}},\endaligned\end{eqnarray}
\begin{eqnarray}\aligned &\label{Th-1-8} N(p)=N_{0}\sum^\infty_{n=0}\frac{(\gamma)_{n\tau,k}}{\Gamma_k(n \alpha+\beta)}\left\{p^{-(n+1)}\sum _{r=0}^{\infty }\left[-\left(\frac{p}{d}\right)^{-\nu}\right]^{r}\right\}.\endaligned \end{eqnarray}
Taking the Laplace inverse of \eqref{Th-1-8}, and using
\begin{eqnarray}\aligned &\label{Th-1-9} L^{-1}\left\{p^{-\nu};t\right\}=\frac{t^{\nu-1}}{\Gamma(\nu)},(R(\nu)>0)\endaligned\end{eqnarray}
we have
\begin{eqnarray}\aligned &\label{Th-1-10} L^{-1}\left\{N(p)\right\}=N_{0}\sum^\infty_{n=0}\frac{(\gamma)_{n\tau,k}}{\Gamma_k(n \alpha+\beta)}L^{-1}\left\{\sum _{r=0}^{\infty }(-1)^{r}d^{\nu r}p^{-(n+1+\nu r)}\right\},\endaligned \end{eqnarray}
\begin{eqnarray}\aligned &\label{Th-1-11} N(t)=N_{0}\sum^\infty_{n=0}\frac{(\gamma)_{n\tau,k}}{\Gamma_k(n \alpha+\beta)}\left\{\sum _{r=0}^{\infty }(-1)^{r}d^{\nu r}\frac{t^{(n+\nu r)}}{\Gamma\left(n+\nu r+1\right)}\right\}\endaligned\end{eqnarray}
\begin{eqnarray}\aligned &\label{Th-1-12}=N_{0}\sum^\infty_{n=0}\frac{(\gamma)_{n\tau,k}}{\Gamma_k(n \alpha+\beta)}t^n\left\{\sum _{r=0}^{\infty}(-1)^{r}\frac{d^{\nu r}t^{\nu r}}{\Gamma\left(\nu r+n+1\right)}\right\},\endaligned\end{eqnarray}
\begin{eqnarray}\aligned &\label{Th-1-13} N(t)=N_{0}\sum^\infty_{n=0}\frac{(\gamma)_{n\tau,k}}{\Gamma_k(n \alpha+\beta)}t^n E_{\nu,n+1}(-d^{\nu}t^{\nu}).\endaligned \end{eqnarray}

\begin{theorem}\label{Th-2} If $  d>0, \nu>0,$ $k\in \mathbb{R}, \alpha, \beta, \gamma\in \mathbb{C}; \Re(\alpha)>0, \Re(\beta)>0, \Re{(\gamma)}>0$ and
$ q\in (0,1)\cup\mathbb{N}$, then the solution of the equation
\begin{eqnarray}\aligned &\label{Th-2-1} N(t)=N_{0}E^{\gamma,\tau}_{k,\alpha,\beta}(d^{\nu}t^{\nu})-d^{\nu}{}_{0}D^{-\nu}_{t}N(t)\endaligned\end{eqnarray}
is given by
\begin{eqnarray}\aligned &\label{Th-2-2} N(t)=N_{0}\sum^\infty_{n=0}\frac{(\gamma)_{n\tau,k}}{\Gamma_k(n \alpha+\beta)}(d^{\nu}t^{\nu})^n E_{\nu,\nu n+1}(-d^{\nu}t^{\nu}),\endaligned\end{eqnarray}
where $ E_{\nu,\nu n+1}(.) $ is the generalized Mittag-Leffler function defined in equation \eqref{2.3}.
\end{theorem}

\begin{theorem}\label{Th-3} If $ a>0, d>0, \nu>0,$ $k\in \mathbb{R}, \alpha, \beta, \gamma\in \mathbb{C}; \Re(\alpha)>0, \Re(\beta)>0, \Re{(\gamma)}>0$ and $ q\in (0,1)\cup\mathbb{N}$, then the solution of the equation
\begin{eqnarray}\aligned &\label{Th-3-1} N(t)=N_{0}E^{\gamma,\tau}_{k,\alpha,\beta}(d^{\nu}t^{\nu})-a^{\nu}{}_{0}D^{-\nu}_{t}N(t)
\endaligned\end{eqnarray}
is given by
\begin{eqnarray}\aligned &\label{Th-3-2} N(t)=N_{0}\sum^\infty_{n=0}\frac{(\gamma)_{n\tau,k}}{\Gamma_k(n \alpha+\beta)}(d^{\nu}t^{\nu})^nE_{\nu,\nu n+1}(-a^{\nu}t^{\nu}),\endaligned\end{eqnarray}
where $ E_{\nu,\nu +1}(.) $ is the generalized Mittag-Leffler function defined in equation \eqref{2.3}.
\end{theorem}

Proof: The proof of \eqref{Th-2-2} and \eqref{Th-3-2} is similar to
that given to prove  \eqref{Th-1-2}.

\section{Special Cases}

If we choose  $q=1$, then \eqref{Th-1-2}, \eqref{Th-2-2} and
\eqref{Th-3-2} reduces to the following:

\begin{corollary} If $ d>0, \nu>0,$ $k\in \mathbb{R}, \alpha, \beta, \gamma\in \mathbb{C}; \Re(\alpha)>0, \Re(\beta)>0$ and  $\Re{(\gamma)}>0$, then the solution of the equation
\begin{eqnarray}\aligned & N(t)-N_{0}E^{\gamma}_{k,\alpha,\beta}(t)=-d^{\nu}{}_{0}D^{-\nu}_{t}N(t)\endaligned\end{eqnarray}
is given by
\begin{eqnarray}\aligned & N(t)=N_{0}\sum^\infty_{n=0}\frac{(\gamma)_{n,k}}{\Gamma_k(n \alpha+\beta)}t^n E_{\nu,n+1}(-d^{\nu}t^{\nu}),\endaligned \end{eqnarray}
where  $ E_{\alpha,\beta}(x) $ is the Mittag-Leffler function defined in equation \eqref{2.3}.

\end{corollary}

\begin{corollary} If $  d>0, \nu>0,$ $k\in \mathbb{R}, \alpha, \beta, \gamma\in \mathbb{C}; \Re(\alpha)>0, \Re(\beta)>0$
and  $\Re{(\gamma)}>0$, then the solution of the equation
\begin{eqnarray}\aligned & N(t)=N_{0}E^{\gamma}_{k,\alpha,\beta}(d^{\nu}t^{\nu})-d^{\nu}{}_{0}D^{-\nu}_{t}N(t)\endaligned\end{eqnarray}
is given by
\begin{eqnarray}\aligned & N(t)=N_{0}\sum^\infty_{n=0}\frac{(\gamma)_{n,k}}{\Gamma_k(n \alpha+\beta)}(d^{\nu}t^{\nu})^nE_{\nu,\nu n+1}(-d^{\nu}t^{\nu}),\endaligned\end{eqnarray}
where $ E_{\nu,\nu n+1}(.) $ is the Mittag-Leffler function defined in equation \eqref{2.3}.
\end{corollary}

\begin{corollary} If $ a>0, d>0, \nu>0,$ $k\in \mathbb{R}, \alpha, \beta, \gamma\in \mathbb{C}; \Re(\alpha)>0, \Re(\beta)>0$ and  $\Re{(\gamma)}>0$, then the solution of the equation
\begin{eqnarray}\aligned &  N(t)=N_{0}E^{\gamma}_{k,\alpha,\beta}(d^{\nu}t^{\nu})-a^{\nu}{}_{0}D^{-\nu}_{t}N(t)
\endaligned\end{eqnarray}
is given by
\begin{eqnarray}\aligned &  N(t)=N_{0}\sum^\infty_{n=0}\frac{(\gamma)_{n,k}}{\Gamma_k(n \alpha+\beta)}(d^{\nu}t^{\nu})^nE_{\nu,\nu n+1}(-a^{\nu}t^{\nu}),\endaligned\end{eqnarray}
where $ E_{\nu,\nu n+1}(.) $ is the generalized Mittag-Leffler function defined in equation \eqref{2.3}.
\end{corollary}


\vskip 3mm \noindent  If we choose  $k=1$, then \eqref{Th-1-2},
\eqref{Th-2-2} and \eqref{Th-3-2} reduces to the following:

\begin{corollary} If $ d>0, \nu>0,$ $\alpha, \beta, \gamma\in \mathbb{C}; \Re(\alpha)>0, \Re(\beta)>0, \Re{(\gamma)}>0$ and $ q\in (0,1)\cup\mathbb{N}$, then the solution of the equation
\begin{eqnarray}\aligned & N(t)-N_{0}E^{\gamma,q}_{\alpha,\beta}(t)=-d^{\nu}{}_{0}D^{-\nu}_{t}N(t)\endaligned\end{eqnarray}
is given by
\begin{eqnarray}\aligned & N(t)=N_{0}\sum^\infty_{n=0}\frac{(\gamma)_{nq,k}}{\Gamma_k(n \alpha+\beta)}t^n E_{\nu,n+1}(-d^{\nu}t^{\nu}),\endaligned \end{eqnarray}
where $ E_{\nu,n+1}(.) $ is the generalized Mittag-Leffler function defined in equation \eqref{2.3}.
\end{corollary}

\begin{corollary} If $ d>0, \nu>0,$ $\alpha, \beta, \gamma\in \mathbb{C}; \Re(\alpha)>0, \Re(\beta)>0, \Re{(\gamma)}>0$ and $ q\in (0,1)\cup\mathbb{N}$, then the solution of the equation
\begin{eqnarray}\aligned &  N(t)-N_{0}E^{\gamma,q}_{\alpha,\beta}(d^{\nu}t^{\nu})=-d^{\nu}{}_{0}D^{-\nu}_{t}N(t)\endaligned\end{eqnarray}
is given by
\begin{eqnarray}\aligned &  N(t)=N_{0}\sum^\infty_{n=0}\frac{(\gamma)_{nq,k}}{\Gamma_k(n \alpha+\beta)}(d^{\nu}t^{\nu})^n E_{\nu,n+1}(-d^{\nu}t^{\nu}),\endaligned \end{eqnarray}
where $ E_{\nu,\nu n+1}(.) $ is the generalized Mittag-Leffler function defined in equation \eqref{2.3}.

\end{corollary}

\begin{corollary} If $ a>0, d>0, \nu>0,$ $\alpha, \beta, \gamma\in \mathbb{C}; \Re(\alpha)>0, \Re(\beta)>0, \Re{(\gamma)}>0$ and $ q\in (0,1)\cup\mathbb{N}$, then the solution of the equation
\begin{eqnarray}\aligned &  N(t)-N_{0}E^{\gamma,q}_{\alpha,\beta}(d^{\nu}t^{\nu})=-a^{\nu}{}_{0}D^{-\nu}_{t}N(t)\endaligned\end{eqnarray}
is given by
\begin{eqnarray}\aligned &  N(t)=N_{0}\sum^\infty_{n=0}\frac{(\gamma)_{nq,k}}{\Gamma_k(n \alpha+\beta)}(d^{\nu}t^{\nu})^nE_{\nu,n+1}(-a^{\nu}t^{\nu}),\endaligned \end{eqnarray}
where $ E_{\nu,\nu n+1}(.) $ is the generalized Mittag-Leffler function defined in equation \eqref{2.3}.
\end{corollary}

\noindent If we choose  $q=1$ and $k=1$, then \eqref{Th-1-2},
\eqref{Th-2-2} and \eqref{Th-3-2} reduces to the following:

\begin{corollary} If $ d>0, \nu>0,$ $\alpha, \beta, \gamma\in \mathbb{C}; \Re(\alpha)>0, \Re(\beta)>0$ and$ \Re{(\gamma)}>0$, then the solution of the equation
\begin{eqnarray}\aligned & N(t)-N_{0}E^{\gamma}_{\alpha,\beta}(t)=-d^{\nu}{}_{0}D^{-\nu}_{t}N(t)\endaligned\end{eqnarray}
is given by
\begin{eqnarray}\aligned & N(t)=N_{0}\sum^\infty_{n=0}\frac{(\gamma)_{n}}{\Gamma(n \alpha+\beta)}t^n E_{\nu,n+1}(-d^{\nu}t^{\nu}),\endaligned \end{eqnarray}
where $ E_{\nu,n+1}(.) $ is the generalized Mittag-Leffler function defined in equation \eqref{2.3}.
\end{corollary}

\begin{corollary} If $  d>0, \nu>0,$ $\alpha, \beta, \gamma\in \mathbb{C}; \Re(\alpha)>0, \Re(\beta)>0$ and$ \Re{(\gamma)}>0$, then the solution of the equation
\begin{eqnarray}\aligned & N(t)=N_{0}E^{\gamma}_{\alpha,\beta}(d^{\nu}t^{\nu})-d^{\nu}{}_{0}D^{-\nu}_{t}N(t)\endaligned\end{eqnarray}
is given by
\begin{eqnarray}\aligned & N(t)=N_{0}\sum^\infty_{n=0}\frac{(\gamma)_{n}}{\Gamma(n \alpha+\beta)}(d^{\nu}t^{\nu})^nE_{\nu,\nu n+1}(-d^{\nu}t^{\nu}),\endaligned\end{eqnarray}
where $ E_{\nu,\nu n+1}(.) $ is the generalized Mittag-Leffler
function defined in equation \eqref{2.3}.
\end{corollary}

\begin{corollary} If $ a>0, d>0, \nu>0,$ $\alpha, \beta, \gamma\in \mathbb{C}; \Re(\alpha)>0, \Re(\beta)>0$ and$ \Re{(\gamma)}>0$, then the solution of the equation
\begin{eqnarray}\aligned &  N(t)=N_{0}E^{\gamma}_{\alpha,\beta}(d^{\nu}t^{\nu})-a^{\nu}{}_{0}D^{-\nu}_{t}N(t)\endaligned\end{eqnarray}
is given by
\begin{eqnarray}\aligned &  N(t)=N_{0}\sum^\infty_{n=0}\frac{(\gamma)_{n}}{\Gamma(n \alpha+\beta)}(d^{\nu}t^{\nu})^nE_{\nu,\nu n+1}(-a^{\nu}t^{\nu}),\endaligned\end{eqnarray}
where $ E_{\nu,\nu n+1}(.) $ is the generalized Mittag-Leffler
function defined in equation \eqref{2.3}.
\end{corollary}

\noindent  If we choose  $q=1, k=1$ and $\gamma=1$, then
\eqref{Th-1-2}, \eqref{Th-2-2} and \eqref{Th-3-2} reduces to the
following:

\begin{corollary} If $ d>0, \nu>0;\alpha, \beta\in \mathbb{C}; \Re(\alpha)>0, \Re(\beta)>0$, then the solution of the equation
\begin{eqnarray}\aligned & N(t)-N_{0}E_{\alpha,\beta}(t)=-d^{\nu}{}_{0}D^{-\nu}_{t}N(t)\endaligned\end{eqnarray}
is given by
\begin{eqnarray}\aligned & N(t)=N_{0}\sum^\infty_{n=0}\frac{t^n}{\Gamma(n \alpha+\beta)} E_{\nu,n+1}(-d^{\nu}t^{\nu}),\endaligned \end{eqnarray}
where $ E_{\nu,n+1}(.) $ is the generalized Mittag-Leffler function defined in equation \eqref{2.3}.
\end{corollary}

\begin{corollary} If $  d>0, \nu>0;\alpha, \beta\in \mathbb{C}; \Re(\alpha)>0, \Re(\beta)>0$, then the solution of the equation
\begin{eqnarray}\aligned &  N(t)=N_{0}E_{\alpha,\beta}(d^{\nu}t^{\nu})-d^{\nu}{}_{0}D^{-\nu}_{t}N(t)\endaligned\end{eqnarray}
is given by
\begin{eqnarray}\aligned &  N(t)=N_{0}\sum^\infty_{n=0}\frac{(d^{\nu}t^{\nu})^n}{\Gamma(n \alpha+\beta)}E_{\nu,\nu n+1}(-d^{\nu}t^{\nu}),\endaligned\end{eqnarray}
where $ E_{\nu,\nu n+1}(.) $ is the generalized Mittag-Leffler
function defined in equation \eqref{2.3}.
\end{corollary}

\begin{corollary} If $ a>0, d>0, \nu>0;\alpha, \beta\in \mathbb{C}; \Re(\alpha)>0, \Re(\beta)>0$, then the solution of the equation
\begin{eqnarray}\aligned &  N(t)=N_{0}E_{\alpha,\beta}(d^{\nu}t^{\nu})-a^{\nu}{}_{0}D^{-\nu}_{t}N(t)\endaligned\end{eqnarray}
is given by
\begin{eqnarray}\aligned &  N(t)=N_{0}\sum^\infty_{n=0}\frac{(d^{\nu}t^{\nu})^n}{\Gamma(n \alpha+\beta)}E_{\nu,\nu n+1}(-a^{\nu}t^{\nu}),\endaligned\end{eqnarray}
where $ E_{\nu,\nu n+1}(.) $ is the generalized Mittag-Leffler
function defined in equation \eqref{2.3}.
\end{corollary}

\noindent If we choose  $q=1, k=1, \gamma=1$ and $\beta=1$, then
\eqref{Th-1-2}, \eqref{Th-2-2} and \eqref{Th-3-2} reduces to the
following:


\begin{corollary} If $ d>0, \nu>0;\alpha\in \mathbb{C}; \Re(\alpha)>0$, then the solution of the equation
\begin{eqnarray}\aligned &  N(t)-N_{0}E_{\alpha}(t)=-d^{\nu}{}_{0}D^{-\nu}_{t}N(t)\endaligned\end{eqnarray}
is given by
\begin{eqnarray}\aligned &  N(t)=N_{0}\sum^\infty_{n=0}\frac{t^n}{\Gamma(n \alpha+1)} E_{\nu,n+1}(-d^{\nu}t^{\nu}),\endaligned \end{eqnarray}
where $ E_{\nu,n+1}(.) $ is the generalized Mittag-Leffler function defined in equation \eqref{2.3}.
\end{corollary}

\begin{corollary} If $  d>0, \nu>0;\alpha\in \mathbb{C}; \Re(\alpha)>0$, then the solution of the equation
\begin{eqnarray}\aligned &   N(t)=N_{0}E_{\alpha}(d^{\nu}t^{\nu})-d^{\nu}{}_{0}D^{-\nu}_{t}N(t)\endaligned\end{eqnarray}
is given by
\begin{eqnarray}\aligned &
N(t)=N_{0}\sum^\infty_{n=0}\frac{(d^{\nu}t^{\nu})^n}{\Gamma(n
\alpha+1)}E_{\nu,\nu n+1}(-d^{\nu}t^{\nu}),\endaligned\end{eqnarray}
where $ E_{\nu,\nu n+1}(.) $ is the generalized Mittag-Leffler
function defined in equation \eqref{2.3}.
\end{corollary}

\begin{corollary} If $ a>0, d>0, \nu>0;\alpha\in \mathbb{C}; \Re(\alpha)>0$, then the solution of the equation
\begin{eqnarray}\aligned &   N(t)=N_{0}E_{\alpha}(d^{\nu}t^{\nu})-a^{\nu}{}_{0}D^{-\nu}_{t}N(t)\endaligned\end{eqnarray}
is given by
\begin{eqnarray}\aligned &   N(t)=N_{0}\sum^\infty_{n=0}\frac{(d^{\nu}t^{\nu})^n}{\Gamma(n \alpha+1)}E_{\nu,\nu n+1}(-a^{\nu}t^{\nu}),\endaligned\end{eqnarray}
where $ E_{\nu,\nu n+1}(.) $ is the generalized Mittag-Leffler
function defined in equation \eqref{2.3}.
\end{corollary}

\noindent  If we choose  $q=1, k=1, \gamma=1, \alpha=0$ and
$\beta=1$, then \eqref{Th-1-2}, \eqref{Th-2-2} and \eqref{Th-3-2}
reduces to the following:


\begin{corollary} If $ d>0, \nu>0$ then the solution of the equation
\begin{eqnarray}\aligned & N(t)-N_{0}e^t=-d^{\nu}{}_{0}D^{-\nu}_{t}N(t)\endaligned\end{eqnarray}
is given by
\begin{eqnarray}\aligned &  N(t)=N_{0}e^t E_{\nu,n+1}(-d^{\nu}t^{\nu}),\endaligned \end{eqnarray}
where $ E_{\nu,n+1}(.) $ is the generalized Mittag-Leffler function defined in equation \eqref{2.3}.
\end{corollary}

\begin{corollary} If $  d>0, \nu>0$ then the solution of the equation
\begin{eqnarray}\aligned &  N(t)=N_{0}e^{(d^{\nu}t^{\nu})}-d^{\nu}{}_{0}D^{-\nu}_{t}N(t)\endaligned\end{eqnarray}
is given by
\begin{eqnarray}\aligned &  N(t)=N_{0}e^{(d^{\nu}t^{\nu})^n} E_{\nu,\nu n+1}(-d^{\nu}t^{\nu}),\endaligned\end{eqnarray}
where $ E_{\nu,\nu n+1}(.) $ is the generalized Mittag-Leffler
function defined in equation \eqref{2.3}.
\end{corollary}

\begin{corollary} If $ a>0, d>0, \nu>0$ then the solution of the equation
\begin{eqnarray}\aligned &  N(t)=N_{0}e^{(d^{\nu}t^{\nu})}-a^{\nu}{}_{0}D^{-\nu}_{t}N(t)\endaligned\end{eqnarray}
is given by
\begin{eqnarray}\aligned &  N(t)=N_{0}e^{(d^{\nu}t^{\nu})^n }E_{\nu,\nu n+1}(-a^{\nu}t^{\nu}),\endaligned\end{eqnarray}
where $ E_{\nu,\nu n+1}(.) $ is the generalized Mittag-Leffler
function defined in equation \eqref{2.3}.
\end{corollary}

\section{Database and Graphical Interpretation}

In this section we establish database for numerical solutions of the
kinetic equations \eqref{Th-1-2}, \eqref{Th-2-2} and \eqref{Th-3-2}
for particular values of the parameters and their graphs and
Mesh-plot are plotted in figure 1 and 2 respectively. Here  we
denote the solutions of equations \eqref{Th-1-2}  as
$N(t)=N(N_0,\gamma,\tau,k,\alpha,\beta,d,\nu,t)$; where
$N_0=.05,\gamma=2, \tau=1, k=2, \alpha=6, \beta=7, d=3, \nu=1$;
those for equation \eqref{Th-2-2}  as
$N(t)=N(N_0,\gamma,\tau,k,\alpha,\beta,d,\nu,t)$; where
$N_0=.05,\gamma=2, \tau=1, k=2, \alpha=6, \beta=7, d=3, \nu=5$ and
those for equation \eqref{Th-3-2}  as
$N(t)=N(N_0,\gamma,\tau,k,\alpha,\beta,d,a,\nu,t)$; where $N_0=.05,
\gamma=2, \tau=1, k=2, \alpha=6, \beta=7, d=3, a=3, \nu=7$. We
obtained  the  database (see Table 1) by using these values. On the
basis of the same values assigned to the parameters in equations
\eqref{Th-1-2}, \eqref{Th-2-2} and \eqref{Th-3-2}, we have the
graphs in figure 1 and 2.
\newline

\end{document}